\documentclass[12pt]{article}

\usepackage{amsfonts}
\usepackage{float}
\usepackage{flafter}

\leftmargin\leftmargini
\labelwidth\leftmargini\advance\labelwidth-\labelsep

\newcounter{defcounter}[section]
{\vspace{0.1cm}\begin{sloppypar}\noindent\stepcounter{defcounter}{\bfseries Definition
      \thesection.\thedefcounter}}%
{\end{sloppypar}\vspace{0.1cm}}
\newtheorem{corollary}{Corollary}[section]
\newtheorem{lemma}{Lemma}[section]
\newtheorem{theorem}{Theorem}[section]
\newtheorem{proposition}{Proposition}[section]

\newcommand{\proof}{{\noindent\bf Proof. }}

\newcommand{\qed}{\hfill $\square$\\~\\}

\begin{document}
\thispagestyle{empty}
\begin{center}
{\Large {\bf Dimension Theory of Linear Solenoids}}
\end{center}
\begin{center}
J. Neunh\"auserer\\~\\
Wachtelpforte. 30, \\ 38640 Goslar, Germany\\
neunchen@aol.com
\end{center}
\begin{abstract}
We develop the dimension theory for a class of linear solenoids,
which have a "fractal" attractor. We will find the dimension of
the attractor, proof formulas for the dimension of ergodic
measures on this attractor and discuss the question whether there
exists a measure of full dimension.
  ~~ \\
{\bf MSC 2000: 37C45, 28D20, 28A80}
\end{abstract}
\section{Introduction}
We consider in this article a class of dynamical systems given by
piece-wise lineare maps acting on a cube. These dynamical system
are very similar to the classical Smale-Williams Solenoids (see
\cite{[KH]}), having a one dimensional unstable and a two
dimensional stable manifold. As the solenoid the systems we study
have a global attractor, which has a complicated "fractal"
geometry. Thus we will discuss here the "fractal" dimension of
this attractor. In Theorem 3.1 we determine the Hausdorff and
Box-counting dimension, using results on the dimension of
self-affine sets in the plane found in \cite{[NE1]}. In the
following we will apply symbolic dynamics and the theory of
hyperbolic dynamics to the class of lineare solenoids. In section
four we will find a coding of the dynamics through a shift on two
symbols on a set of full measure. This allows us to to find a
representation of all ergodic measures for the systems as images
of shift ergodic measures under the coding map. In section five we
will demonstrate the existence of Lyapunov exponents and Lyapunov
charts for lineare solenoids with respect to any ergodic measure.
This is the background we need to apply the general dimension
theory of hyperbolic measures, see \cite{[LY]} and \cite{[BPS]}.
Using this theory we will show that ergodic measure for our
systems are exact dimensional. Moreover we will find a formula for
the dimension of ergodic measures in terms of entropy and Lyapunov
exponents and the dimension of transversal measures (see Theorem
6.1). For Bernoulli measures this formula yield an explicit
expression (see Corollary 6.2) for the dimension in terms of
self-similar measure studied in \cite{[NE1]} . In the last section
of this article we will discuss the question whether there exists
an ergodic measure of full dimension, which means that the
dimension of the ergodic measure equals the dimension of the
attractor. This question is widely open in the dimension theory of
dynamical systems. It is of particular interest since ergodic
measures of full dimension are of great geometrical significants,
describing the long term behavior of orbits on the whole attractor
in the dimensional theoretical sense . Results of Manning and
McClusky \cite{[MM]} show that in the case of horseshoes
diffeomorphisms there does not exist an ergodic measure of full
dimension in general. One can not maximize the dimension in the
stable and in the unstable direction at the same time. In
\cite{[NE2]} we demonstrate that for generalized Baker`s
transformations there exists parameter domains for which a measure
of full dimension exists and parameter domains where the dimension
of the invariant set can not even be approximated by the dimension
of ergodic measures. We observe the same phaenomenon in the case
of lineare solenoids. We will show that there are manifolds in the
parameter domain where there is a measure of full dimension and
manifold where the variational principle of dimension does not
hold (see Theorem 7.1). At the end of this paper the reader, who
is not familiar
with dimension theory, will find an appendix containing a short introduction to this field.\\~\\
{\bf Acknowledgment}\newline\newline I wish to thank the
supervisor of my PhD J\"org Schmeling who helped me a lot to find
the results presented here.
\section{Linear solenoids}
Let $\mathbb{W}=[-1,1]^{3}$. We consider the following class of
piecewise affine maps $f_{v}:\mathbb{W}\longmapsto \mathbb{W}$
given by
\begin{displaymath}
{ f_{v} (x,y,z) } = {\lbrace {\begin{array}{cc}
 (2x-1,\beta_{1} y +(1-\beta_{1}),\tau_{1} z+(1-\tau_{1}))\quad \mbox{if} \quad x \geq 0   \\
 (2x+1,\beta_{2} y -(1-\beta_{2}),\tau_{2} z-(1-\tau_{2}))\quad \mbox{if} \quad x<0
\end{array} }}
\end{displaymath}
where we assume
\[ v=(\beta_{1},\beta_{2},\tau_{1},\tau_{2})\in(0,1)^{4}\quad\mbox{ and }\quad \tau_{1}+\tau_{2}< 1. \]
\begin{center}
\unitlength=1.00mm \special{em:linewidth 0.4pt}
\linethickness{0.4pt}

\end{center}
{\bf Figure 1:} The action of $f_{v}$ on the cube $\mathbb{W}$.
\\~\\~Obvoiusly the maps $f_{v}$ are invertible and there is a global
Attraktor for the maps given by
\[
\Lambda_{v}=\mbox{closure}(\bigcap_{k=0}^{\infty}f^{n}_{v}(\mathbb{W})).
\]
We call the system $(\Lambda_{v},f_{v})$ a linear Solenoid. We see
that this system is quite similar to the classical Smale-Williams
Solenoid $(\Delta_{\beta,\tau},g_{\beta,\tau})$, see \cite{[KH]}.
The Smale-Williams Solenoid is constructed by a family of maps
$g_{\beta,\tau}:\mathbb{T}^2\longmapsto \mathbb{T}^2$ on the full
torus $\mathbb{T}^2= \mathbb{S}^{1}\times\mathbb{D}^2$ defined by
\[g(\phi,y,z)=(2\phi \mbox{ mod }2\pi,\beta
y+1/2\cos(2\pi\phi),\tau z+1/2\sin(2\pi\phi)) \]
\[ \mbox{with} \qquad \beta,\tau\in (0,1/2). \]
$g_{\beta,\tau}$ has the global attractor
\[ \Delta_{\beta,\tau}=\bigcap_{n=0}^{\infty}g_{\beta,\tau}^{n}(\mathbb{T}^{2}). \]
\begin{figure}[H]
\end{figure}
{\bf Figure 2} The action of $g_{\beta,\tau}$ on the full Torus
$\mathbb{T}^{2}$.~~\\~\\ In fact the systems
$(\Delta_{\beta,\tau},g_{\beta,\tau})$ and $(\Lambda_{v},f_{v})$
have similar properties. Both systems are expanding in the first
coordinate direction with expansion rate $\log 2$ and contracting
in the two other coordinate direction. Moreover both maps are
invertible and into with a global attractor. The classical
Solenoid is hyperbolic and conjugated to the full Shift on two
symbols, see \cite{[KH]}. In section four we will show that our
linear Solenoids are up to a set of measure zero as well
hyperbolic and conjugated to the full shift on two symbols. The
last similarity is that both $\Lambda_{v}$ and
$\Delta_{\beta,\tau}$ have a complicated non smooth geometry.
Dimensional theoretical properties of the classical Solenoid where
extensity studied, see \cite{[BO]},\cite{[SI]} or \cite{[PE]}. We
will develop here the dimensiontheory for lineare solenoids. In
the next chapter we will present our results on the dimension of
the attractor $\Lambda_{v}$.
\section{Dimension of the attractor}
We first give here an simple description of the attractor
$\Lambda_{v}$ using iterated function systems, see \cite{[FA]}.
\begin{proposition}
We have
\[ \Lambda_{v}=[-1,1]\times \Lambda_{v}^{s} \]
where $\Lambda_{v}^{s}$ is the unique compact set fulfilling
\[ \Lambda_{v}^{s}=T_{\beta_{1},\tau_{1}}(\Lambda_{v}^{s})\cup
T_{\beta_{2},\tau_{2}}(\Lambda_{v}^{s}) \] with
$T_{\beta_{1},\tau_{1}},T_{\beta_{2},\tau_{2}}:[-1,1]^{2}\longmapsto
[-1,1]^{2}$ given by
\[ T_{\beta_{1},\tau_{1}}(y,z)=(\beta_{1} y +(1-\beta_{1}),\tau_{1} z+(1-\tau_{1}))\]
\[ T_{\beta_{2},\tau_{2}}(y,z)=(\beta_{2} y +(1-\beta_{2}),\tau_{2} z+(1-\tau_{2})).\]
\end{proposition}
\proof Let $T_{1}:=T_{\beta_{1},\tau_{1}}$ and
$T_{2}=T_{\beta_{2},\tau_{2}}$. We have
\[ \mbox{closure}(f_{v}(\mathbb{W}))=[-1,1]\times T_{1}(\mathbb{W})\cup
[-1,1]\times T_{2}(\mathbb{W})\] and hence
\[ \mbox{closure}(f_{v}^{n}(\mathbb{W})=[-1,1] \times
\bigcup_{s_{1},\dots,s_{n}\in\{1,2\}}T_{s_{1}}\circ T_{s_{2}}\circ
\dots \circ T_{s_{2}}(\mathbb{W}) \] Now let
\[ \Lambda^{s}_{v}=\bigcap_{n=1}^{\infty}\bigcup_{s_{1},\dots,s_{n}\in\{1,2\}}T_{s_{1}}\circ T_{s_{2}}\circ
\dots \circ T_{s_{2}}(\mathbb{W})\] By this definition we get
\[ \Lambda_{v}=[-1,1]\times \Lambda_{v}^{s}. \]
Moreover $\Lambda_{v}^{s}$ is compact with
\[
\Lambda_{v}^{s}=T_{1}(\Lambda_{v}^{s})\cup T_{2}(\Lambda_{v}^{s})
\]
Uniqueness of $\Lambda_{v}^{s}$ with this property follows from
\cite{[HU]}. \qed Our results on the dimension of the attractor
$\Lambda_{v}$ is now mainly a consequence of our results on the
self-affine sets $\Lambda_{v}^{s}$ given in \cite{[NE1]} and and
\cite{[PE]} . In the following we denote by $\dim_{B}A$ the
box-Counting dimension and by $\dim_{H}A$ the Hausdorff dimension
of a set $A$; we refer to the appendix of this work for the
definition of these quantities.
\begin{theorem}
Let $v=(\beta_{1},\beta_{2},\tau_{1},\tau_{2})\in(0,1)^{4}$ with
$\beta_{1}+\beta_{2}>\tau_{1}+\tau_{2}$.\\~\\
If $\beta_{1}+\beta_{2}<1$ we have
\[ \dim_{B}\Lambda_{v}=\dim_{H}\Lambda_{v}=d+1\]
where $d$ is the solution of
\[ \beta_{1}^{d}+\beta_{2}^{d}=1. \]
If $\beta_{1}+\beta_{2}\ge 1$ we have
\[ \dim_{B}\Lambda_{v}=d+2\]
where $d$ is the solution of
\[ \beta_{1}\tau_{1}^{d}+\beta_{2}\tau_{2}^{d}=1. \]
Moreover for almost all $\beta_{1},\beta_{2}<0.649$ we have
\[ \dim_{H}\Lambda_{v}=\dim_{B}\Lambda_{v}.\]
\end{theorem}
\proof By proposition 8.1 of the appendix and proposition 3.1 we
have
\[ \dim_{H/B}\Lambda_{v}=\dim_{H/B}\Lambda_{v}^{s}+1. \]
If $\beta_{1}+\beta_{2}<1$ we have by example 16.3 of \cite{[PE]} $\dim_{B}\Lambda_{v}=\dim_{H}\Lambda_{v}=d$. \\~\\
If $\beta_{1}+\beta_{2}\ge 1$ we get by theorem III of
\cite{[NE1]} $\dim_{B}\Lambda^{s}_{v}=d+2$ and generically
$\dim_{H}\Lambda^{s}_{v}=\dim_{B}\Lambda^{s}_{v}$ under the
assumption that $\beta_{1},\beta_{2}<0.649$. \qed The condition
$\beta_{1},\beta_{2}<0.649$ in the last theorem is due to the
technique we used in \cite{[NE1]}. We do not belief that this
condition is essential, also we were not able to omit it. The
identity of Box-Counting and Hausdorff dimension in the last
statement does not hold in general. In \cite{[NE2]} we described
numbertheoretical exceptions in the symmetric case
$\beta_{1}=\beta_{2}$.
\section{Shift coding of the Dynamics}
We need some notation to introduce a symbolic coding of the
dynamics of system $(\Lambda_{v},f_{v})$.~\\~\\
Let $\Sigma=\{-1,1\}^{\mathbb{Z}}$ be the Shift space. With the
product metric defined by
\[ d(s,t)=\sum_{k=-\infty}^{\infty}
|s_{k}-t_{k}|2^{-|k|}
\]
$\Sigma$ becomes a perfect, totally disconnected and compact
metric space; see \cite{[DGS]}.  The forward shift map $\sigma$ on
$\Sigma$ is given by $\sigma((s_{k}))=(s_{k+1})$, the backward
shift $\sigma^{-1}$ is given by $\sigma((s_{k}))=(s_{k-1})$.
\newline\newline
For a sequence $s\in\Sigma$ and $\gamma_{1},\gamma_{2}\in(0,1)$ we
define a map \[ \hat\pi_{\gamma_{1},\gamma_{2}}:\Sigma^{+}
\longrightarrow [\frac{-\gamma_{2}}{1-\gamma_{2}}
\frac{\gamma_{1}}{1-\gamma_{1}}]\] by
\[  \hat\pi_{\gamma_{1},\gamma_{2}}(s)=\sum_{k=0}^{\infty}s_{k}\gamma_{2}^{\sharp(s,k)}\gamma_{1}^{\bar\sharp(s,k)}
\]
where
\[ \sharp(s,k)=\mbox{Cardinality}\{s_{i}|s_{i}=-1~~i=1,\dots,k\}    \]
\[ \bar\sharp(s,k)=\mbox{Cardinality}\{s_{i}|s_{i}=+1~~i=1,\dots,k\} .   \]
Let $L_{\gamma_{1},\gamma_{2}}$ be the monoton increasing linear
map form $[\frac{-\gamma_{2}}{1-\gamma_{2}}
\frac{\gamma_{1}}{1-\gamma_{1}}]$ onto $[-1,1]$ and let
$\pi_{\gamma_{1},\gamma_{2}}=L_{\gamma_{1},\gamma_{2}}\circ
\hat\pi_{\gamma_{1},\gamma_{2}} $. Moreover we define the map of
the signed dyadic expansion
\[ i:\Sigma\longmapsto [-1,1]\]
by
\[ i(s)=\sum_{k=1}^{\infty}s_{-k}(1/2)^{k}. \]
For $v=(\beta_{1},\beta_{2},\tau_{1},\tau_{2})\in (0,1)^{4}$ we
define the coding map
\[ \pi_{v}:\Sigma\longmapsto \Lambda_{v} \]
by
\[ \pi_{v}(s)=(i(s),\pi_{\beta_{1},\beta_{2}}(s),\pi_{\tau_{1},\tau_{2}}(s)). \]
By this definitions we obviously have:
\begin{proposition}
$\pi_{v}$ is continuous and onto $\Lambda_{v}$. Moreover the map
is bijective from
\[\bar\Sigma:=\Sigma\backslash\{(s_{k})|
\exists k_{0} \forall k \le k_{0}\in\mathbb{Z} :s_{k}=1
\rbrace\cup\{(1)\}\] onto $\Lambda_{v}$ and we have
\[  \forall s\in \bar\Sigma~~:~~ \pi_{v}(\sigma^{-1}(s))=f_{v}(\pi_{v}(s)). \]
\end{proposition}
We can now represent all ergodic measures of the system
$(\Lambda_{v},f_{v})$ using the coding map $\pi$. Again we need
some notations. Given a compact metric space $X$ we denote the set
of all Borel probability measures on $X$ by $M(X)$. With the
$\mbox{weak}^{*}$ topology $M(X)$ becomes a compact, convex and
metricable space. If $T$ is a Borel measurable transformation on
$X$ we call a measure $\mu$ $T$-invariant if \[ T(\mu):=\mu\circ
T^{-1}=\mu.\] The set of all invariant measures forms  a compact,
convex and nonempty subset of $M(X)$. An invariant measure $\mu$
is called ergodic if
\[ T^{-1}(B):=B\Rightarrow \mu(B)\in\{0,1\} \] hold
for all Borel subsets $B$ of the space $X$. The set of all ergodic
measures
\[ \mathcal{E}(X,T):=\{\mu\in M(X)|\mu~~T\mbox{-ergodic}\}\] is nonempty,
convex and compact with respect to the weak$^{\star}$ topology. It
consists of the extreme points of the set of invariant measures.
By $b^{p}$ for $p\in(0,1)$ we denote the Bernoulli measure on
$\Sigma$, which is the product of the discrete measure giving $1$
the probability $p$ and $-1$ the probability $(1-p)$. The
Bernoulli measures are ergodic with respect to forward and
backward shifts. Given $b^{p}$ on $\{-1,1\}^{ \mathbb{Z}^{-} }$ we
define the corresponding Bernoulli measure $\ell^{p}$ on $[-1,1]$
by $\ell^{p}=i(b^{p})$. For the basic facts in ergodic theory
mentioned here we refer to \cite{[DGS]}, \cite{[WA]} or
\cite{[KH]}. Proposition 4.1 directly implies:
\begin{proposition}
The map
\[
\mu\longmapsto\mu_{v}:=\mu\circ\pi_{v}^{-1}\] is a affine
homeomorphism from $M(\Sigma,\sigma)$ onto $M(\Lambda_{v},f_{v})$.
Moreover $b^{p}_{v}$ is a product of the Bernoulli measure on
$\Lambda_{v}^{s}$ with $\ell^{p}$.
\end{proposition}
\section{Hyperbolicity}
We will show here that there exists expansion and contraction
rates (Lyapunov exponents) on the Solenoid $(\Lambda_{v},f_{v})$
for a set of full measure with respect to any ergodic measure
$\mu\in \mathcal{E}(\Lambda_{v},f_{v})$.
\begin{lemma} There is a subset
$\Omega_{v}\subseteq\Lambda_{v}$ which has full measure for all
$\mu_{v}\in \mathcal{E}(\Lambda_{v},f_{v})$ such that $f_{v}$ is a
bijection on $\Omega_{v}$ and $f_{v}$ is differentiable for all
${\bf x} \in\Omega_{v}$ with
\[ D_{{\bf  x}} f_{v}=\left( \begin{array}{ccc}
2 &0&0\\
0 & \beta_{1} & 0\\
0& 0 & \tau_{1}
\end{array}\right)
\quad\mbox{if}\quad y>0\quad\mbox{and}\quad
 D_{{\bf x}} f_{v}=\left( \begin{array}{ccc}
2 &0&0\\
0 & \beta_{2} & 0\\
0& 0 & \tau_{2}
\end{array}\right)\quad\mbox{if}\quad y<0.
\]
\end{lemma}
\proof Denote the singularity $\{0\}\times [-1,1]^{2}$ of the map
$f_{v}$  by $S$ and define the set $\Omega_{v}$ by
\[\Omega_{v}=\bigcap_{n=-\infty}^{\infty}f_{v}^{n}(\mathbb{W}\backslash S). \]
By definition we have $ f_{v}(\Omega_{v})=\Omega_{v}$ and since
$f_{\vartheta}$ is injective it is in fact a bijection on
$\Omega_{v}$. Moreover if $(x,y,z)\in \Omega_{v}$ then
$(x,y,z)\not\in S$ and hence $f_{v}$ is differentiable and has
obviously the derivative that we stated in the lemma. It remains
to show that $\mu_{v}(\Omega_{v})=1$. By elemental calculations we
see that
\[\Omega_{v}=(\{(x,y,z)\in \Lambda_{v}|y\not=1,\quad y\not=-1\}\cup\{(1,1,1),(-1,-1,-1)\})\backslash \bigcup_{n=0}^{\infty}f^{-n}(S).\]
Since $\mu_{v}$ is invariant and the union in the expression above
is disjoint it has zero measure. It remains to show that
$\mu_{v}(\{1\}\times[-1,1]\times[-1,1])=\mu_{v}(\{(1,1,1)\}$ and
$\mu_{v}(\{-1\}\times[-1,1]\times[-1,1])=\mu_{v}(\{(-1,-1,-1)\}$.
But this is obvious since $f_{v}$ is just a contraction with fixed
point $(1,1,1)$ resp. $(-1,-1,-1)$ on the sets
$\{1\}\times[-1,1]\times[-1,1]$ resp. $\times\{1\}\times[-1,1]
\times[-1,1]$ . \qed We now define linear subspaces of
$\mathbb{W}$ by
\[ \mathbb{E}^{u}=<\left(\begin{array}{ccc} 1\\0\\0 \end{array}\right)>
\quad \mathbb{E}^{s}=<\left(\begin{array}{ccc} 0\\0\\1
\end{array}\right), \left(\begin{array}{ccc} 0\\1\\0
\end{array}\right)>
\]
\[ \mathbb{E}^{ss}=<\left(\begin{array}{ccc} 0\\0\\1 \end{array}\right)>
\quad \mathbb{E}^{ws}=<\left(\begin{array}{ccc} 0\\1\\0
\end{array}\right)>.
\]
Given a Borel measure $\mu$ on $\Sigma$ and
$\gamma_{1},\gamma_{2}\in(0,1)$ we write
\[ \Xi^{\mu}_{\gamma_{1},\gamma_{2}}=\mu(\{s\in\Sigma|s_{0}=1\})\log\gamma_{1}+\mu(\{s\in\Sigma|s_{0}=-1\})\log\gamma_{2}. \]
\begin{proposition}
Given $\mu\in \mathcal{E}(\Sigma,\sigma)$ we have for
$\mu_{v}$-almost all ${\bf x}\in \Lambda_{v}$
\[
\lim_{n\longrightarrow\infty}\frac{1}{n}\log||D_{{\bf
x}}f^{n}_{\vartheta}\overrightarrow{v}||=\log 2 \qquad \forall
\overrightarrow{v}\in \mathbb{E}^{u}\]
\[\mbox{If }~~ \Xi^{\mu}_{\beta_{1},\beta_{2}}\ge\Xi^{\mu}_{\tau_{1},\tau_{2}}~:~
\lim_{n\longrightarrow\infty}\frac{1}{n}\log||D_{{\bf
x}}f^{n}_{\vartheta}\overrightarrow{v}||=\lbrace
\begin{array}{ccc}
\Xi^{\mu}_{\beta_{1},\beta_{2}} & \mbox{if} & \overrightarrow{v}\in \mathbb{E}^{s}\backslash \mathbb{E}^{ss}\\
\Xi^{\mu}_{\tau_{1},\tau_{2}}   & \mbox{if} & \overrightarrow{v}\in \mathbb{E}^{ss}\\
\end{array}
\]
\[
\mbox{If }~~
\Xi^{\mu}_{\beta_{1},\beta_{2}}\le\Xi^{\mu}_{\tau_{1},\tau_{2}}~:~
\lim_{n\longrightarrow\infty}\frac{1}{n}\log||D_{{\bf
x}}f^{n}_{\vartheta}\overrightarrow{v}||=\lbrace
\begin{array}{ccc}
\Xi^{\mu}_{\tau_{1},\tau_{2}} & \mbox{if} & \overrightarrow{v}\in \mathbb{E}^{s}\backslash \mathbb{E}^{ws}\\
\Xi^{\mu}_{\beta_{1},\beta_{2}}   & \mbox{if} & \overrightarrow{v}\in \mathbb{E}^{ws}\\
\end{array}\]
\end{proposition}
\proof By lemma 5.1 we have for $\mu_{v}$-almost all ${\bf x}\in
\Lambda_{v}$ and all $n>0$
\[ \log||D_{{\bf x}}f^{n}_{v}(\left(\begin{array}{ccc} 0\\y\\0 \end{array}\right))||=n\log2+\log y.\]
This implies our claim about $\mathbb{E}^{u}$. Now consider
$\mathbb{E}^{s}$.  By lemma 5.1, proposition 4.1 and 4.2 we have
for $\mu_{v}$-almost all ${\bf x}\in \Lambda_{v}$ and all $n>0$
\[ \log||D_{{\bf x}}f^{n}_{v}(\left(\begin{array}{ccc} x\\0\\z \end{array}\right))||=\log\sqrt{(x\beta_{1}^{n-\bar\sharp_{n}(s)+1}\beta_{2}^{\bar\sharp_{n}(s)})^{2}+(z\tau_{1}^{n-\bar\sharp_{n}(s)+1}\tau_{2}^{\bar\sharp_{n}(s)})^{2}}\]
where $s=(s_{k})=\pi_{v}^{-1}({\bf x})$ and $\bar\sharp_{n}(s)$
counts the number of entries in the set
$\{s_{0},s_{-1},\dots,s_{-n}\}$ that are $-1$. Applying Birkhoffs
ergodic theorem (see 4.1.2 of \cite{[KH]}) to the functions
\[f_{ws}(s)=\{
\begin{array}{ccc} \log\beta_{1} & \mbox{if} & s_{0}=1\\
                  \log\beta_{2} & \mbox{if} & s_{0}=-1
\end{array}
\qquad f_{ss}(s)=\{
\begin{array}{ccc}\log\tau_{1} & \mbox{if} & s_{0}=1\\
                 \log\tau_{2} & \mbox{if} & s_{0}=-1\end{array}
         \]
now yields the desired result.
\begin{flushright}
$\Box$
\end{flushright}
Proposition 5.1 means that Lyapunov exponents exists almost
everywhere for the ergodic systems $(\Lambda_{v},f_{v},
\mu_{\vartheta})$. $\mathbb{E}^{u}$ is the unstable direction with
Lyapunov exponent $\log 2$ and $\mathbb{E}^{s}$ is the stable
direction with exponent $\Xi^{\mu}_{\beta_{1},\beta_{2}}$ or
$\Xi^{\mu}_{\tau_{1},\tau_{2}}$ depending on which quantity is
bigger. Accordingly $\mathbb{E}^{ss}$ (resp. $E^{ws}$) is the
strong stable direction with Lyapunov exponent
$\Xi^{\mu}_{\tau_{1},\tau_{2}}$ (resp.
$\Xi^{\mu}_{\beta_{1},\beta_{2}}$). In order to guarantee the
existence of Lyapunov charts associated with the Lyapunov
exponents we have to show that the set of points that does not
approach the singularity $S=\{0\}\times[-1,1]\times[-1,1]$ with
exponential rate has full measure, see \cite{[ST1]}. Precisely we
have:
\begin{proposition}
Given $\mu\in \mathcal{E}(\Lambda_{v},f_{v})$ we have for all
$\epsilon>0$
\[ \mu_{v}(\{{\bf x}\in\Lambda_{v}|\exists l>0~\forall n>0~:~d(f^{n}({\bf x}),S)>(1/l)e^{-\epsilon n}\})=1.\]
\end{proposition}
\proof
Fix $\epsilon>0$. First note that it is sufficient if we
show
\[ \mu_{\vartheta}(\{{\bf x}\in\Lambda_{\vartheta}|\exists (n_{k})_{k\in\mathbb{N}}\longrightarrow \infty~\forall k>0~:~d(f^{n_{k}}({\bf x}),S)\le e^{-\epsilon n_{k}}\})=0\]
because if we have for a point ${\bf x}$ that $\exists
n_{0}\forall n>n_{0}$ $d(f^{n}({\bf x}),S)>e^{-\epsilon n}$ then
there exists $l>0$ such that $d(f^{n}({\bf
x}),S)>(1/l)e^{-\epsilon n}$ $\forall n>0$. \newline By
proposition 4.2 and the definition of the measure $\mu_{v}$ this
assertion is equivalent to $\mu(N)=1$ where
\[ N:=\{s\in\hat\Sigma|\exists
(n_{k})_{k\in\mathbb{N}}\longrightarrow \infty~\forall
k>0~d(\sigma^{-n_{k}}(s),\tilde S)\le e^{-\epsilon n_{k}}\}\] and
$\tilde S$ is the singularity in the symbolic coding, i.e.
\[ \tilde S=\{s\in\Sigma|s_{-1}=1 \quad s_{k}=-1~\forall k<-1\}.
\] We will prove this. If $s\in N$ we have
\[ d(\sigma^{-n_{k}}(s),\tilde S)\le e^{-\epsilon n_{k}}\qquad \forall
k>0 \] By the definition of the metric $d$ this implies
\[ (\sigma^{i}(s))_{-2}\not = 1~\mbox{ for }~~i=n_{k},\dots ,n_{k}+\lceil c\epsilon n_{k}\rceil-1 ~~~\forall k>0. \]
where the constant $c$ is  independent of $\epsilon$, $n_{k}$ and
$\underline{s}$. Thus we have:
\[ N\subseteq\{s|\exists (n_{k})_{k\in\mathbb{N}}\longrightarrow \infty~\forall k>0~:~ (\sigma^{i}(s))_{-2}\not = 1~~i=n_{k},\dots ,n_{k}+\lceil c\epsilon n_{k}\rceil-1\}.\]
Applying lemma 7.1. of \cite{[ST2]} for the ergodic system
$(\Sigma,\sigma,\mu)$ we obtain $\mu(N)=0$.\qed
\section{Dimension formulae for ergodic measures}
Our results in the last section demonstrate that we may apply the
general dimension theory for hyperbolic systems to the linear
solenoids $(\Lambda_{v},f_{v},\mu_{v})$. By proposition 5.1 and
proposition 5.2 our systems fall into the class of generalized
hyperbolic attractors in the sense of Schmeling and Troubetzkoy
\cite{[ST1]}. Usually the dimension theory of ergodic measures is
stated in the context of $C^{2}$-diffeomorphisms in order to
guarantee the existence of Lyapunov exponents and charts. But
invertibility and the existence of Lyapunov exponents and charts
almost everywhere is enough to apply this theory. We refer to
section 4 of \cite{[ST1]} for this fact.
\newline First we define here stable partitions $\mathbb{W}^{s}$ and unstable partitions $\mathbb{W}^{u}$
of $\mathbb{W}$ by the partition elements
\[\mathbb{W}^{s}({\bf x})=\{x\}\times[-1,1]\times [-1,1]
\qquad \mathbb{W}^{u}({\bf x})=[-1,1]\times \{y\} \times\{z\}
\] where ${\bf x}=(x,y,z)\in \Lambda_{v}$. Given $\mu_{v}\in
\mathcal{E}(\Lambda_{v},f_{v})$ we have conditional measures
$\mu_{v}^{s}({\bf x})$ on the partition $\mathbb{W}^{s}$ and
conditional measures $\mu_{v}^{u}({\bf x})$ on the partition
$\mathbb{W}^{u}$. These measures are unique $\mu_{v}$-almost
everywhere fulfilling the relations:
\[\mu_{v}(B)=\int\mu^{s}_{v}({\bf x })
(B\cap \mathbb{W}^{s}({\bf x}))d\mu_{v}\quad \mbox{resp.}\quad
\mu_{v}(B)=\int\mu^{u}_{v}({\bf x }) (B\cap \mathbb{W}^{u}({\bf
x}))d\mu_{v}\] for all Borel sets $B$ in $\mathbb{W}$. We refer to
\cite{[RO]} for information about
conditional measures on measurable partitions.\newline\\
To formulate our next theorem let us denote the entropy of an
ergodic measure $\mu$ by $h(\mu)$. We recommend \cite{[WA]} for an
introduction to theory of this invariant. Moreover we denote the
dimension of a measure by $\dim\mu$, so $\mu$ is exact
dimensional. In the end of the appendix the reader finds an
introduction of this quantity.\\~\\
Applying the dimension theory  of hyperbolic systems by Barreira,
Schmeling and Pesin \cite{[BPS]} and Ledrappier Young \cite{[LY]}
to the system $(\Lambda_{v},f_{v},\mu_{v})$ we obtain the
following theorem.
\begin{theorem}
For all $\mu\in \mathcal{E}(\Sigma,\sigma)$ the ergodic measures
$\mu_{v}\in \mathcal{E}(\Lambda_{v},f_{v})$ and the conditional
measures $\mu^{s}_{v}$ and $\mu^{u}_{v}$ are exact dimensional
with
\[ \dim\mu_{v}=\dim \mu^{u}_{v}+\dim \mu^{s}_{v}\]
Moreover we have
\[ \dim\mu^{u}_{v}=h(\mu)/\log 2\]
and
\[\dim\mu^{s}_{v}=\frac{h(\mu)}{-\Xi^{\mu}_{\tau_{1},\tau_{2}}}+(1-\frac{\Xi^{\mu}_{\beta_{1},\beta_{2}}}{\Xi^{\mu}_{\tau_{1},\tau_{2}}})\dim pr_{{\it y}}(\mu_{v})\]
if
$\Xi^{\mu}_{\beta_{1},\beta_{2}}\ge\Xi^{\mu}_{\tau_{1},\tau_{2}}$,
resp.
\[\dim\mu^{s}_{v}=\frac{h(\mu)}{-\Xi^{\mu}_{\beta_{1},\beta_{2}}}+(1-\frac{\Xi^{\mu}_{\tau_{1},\tau_{2}}}{\Xi^{\mu}_{\beta_{1},\beta_{2}}})\dim pr_{{\it z}}(\mu_{v})\]
if
$\Xi^{\mu}_{\beta_{1},\beta_{2}}<\Xi^{\mu}_{\tau_{1},\tau_{2}}$.
Here $pr$ denotes the projection of the measure on second resp.
third coordinate axis.
\end{theorem}
\proof Exact dimensionality of the measures follows directly from
\cite{[DGS]} given our results in section five. The dimension
formula for $\mu^{u}_{v}$ follows directly from theorem $C'$ of
\cite{[LY]}. For the second formula we need an additional
argument. If
$\Xi^{\mu}_{\beta_{1},\beta_{2}}\ge\Xi^{\mu}_{\tau_{1},\tau_{2}}$
\[ \mathbb{W}^{ss}({\bf x})=\{x\}\times\{y\} \times[-1,1]  \]
forms the strong stable partition. We have conditional measures
$\mu_{v}^{ss}({\bf x})$ on $\mathbb{W}^{ss}$. These measures are
unique $\mu_{v}$-almost everywhere fulfilling the relation:
\[\mu_{v}(B)=\int \mu^{ss}_{v}({\bf x })
(B\cap \mathbb{W}^{ss}({\bf x}))d \mu_{v} \] for all Borel sets
$B$ in $\mathbb{W}$. From the uniqueness of the conditional
measures we have for $\mu_{v}$-almost all ${\bf x}$
\[ \mu^{s}_{v}({\bf x})(B)=\int \mu_{v}^{ss}({\bf x})
(B\cap \mathbb{W}^{ss}({\bf x})d pr_{y}(\mu_{v})\] for all Borel
sets B in $W^{s}(\bf x)$. This statement means that the
transversal measures in the sense of  \cite{[LY]} of the nested
partitions $\mathbb{W}^{s}$ and $\mathbb{W}^{ss}$ are given by
$pr_{y} \mu_{v}$. Now the second formula follows again from
theorem $C'$ of \cite{[LY]}. The third formula is proved the same
way just noticing that the strong stable partition is given by
\[ \mathbb{W}^{ss}({\bf x})=\{x\}\times[-1,1] \times \{z\}  \]
in this case. \qed The formula for the conditional measures
$\mu^{s}_{v}$ in theorem 6.1 is known in the dimension theory as
{\it Ledrappier-Young formula}. For Bernoulli measures $b^{p}\in
M(\Sigma,\sigma)$ we get by theorem 6.1 the following explicit
dimension formulas for the measures $b_{v}^{p}\in
M(\Lambda_{v},f_{v})$.
\begin{corollary}
\[ \dim b_{v}^{p}=\frac{p\log
p+(1-p)\log(1-p)}{\log 2}+ \frac{p\log
p+(1-p)\log(1-p)}{p\log\tau_{1}+(1-p)\log\tau_{2}}
 \] \[+(1-\frac{p\log\beta_{1}+(1-p)\log\beta_{2}}{p\log\tau_{1}+(1-p)\log\tau_{2}})\dim \pi_{\beta_{1},\beta_{2}}(b^{p}).\]
if $p\log\beta_{1}+(1-p)\log\beta_{2}\ge
p\log\tau_{1}+(1-p)\log\tau_{2}$
\[ \dim b_{v}^{p}=\frac{p\log
p+(1-p)\log(1-p)}{\log 2}+ \frac{p\log
p+(1-p)\log(1-p)}{p\log\beta_{1}+(1-p)\log\beta_{2}}
 \] \[+(1-\frac{p\log\tau_{1}+(1-p)\log\tau_{2}}{p\log\beta_{1}+(1-p)\log\beta_{2}})\dim \pi_{\tau,\tau_{2}}(b^{p}).\]
if $p\log\beta_{1}+(1-p)\log\beta_{2}<
p\log\tau_{1}+(1-p)\log\tau_{2}$.
\end{corollary}
\proof It is well known in ergodic theory that
\[ h(b^{p})=-(p\log p+(1-p)\log(1-p)) \] see \cite{[DGS]}.
Furthermore we obviously have
\[ \Xi^{b^{p}}_{\gamma_{1},\gamma_{2}}=(p\log\gamma_{1}+(1-p)\log\gamma_{2})\]
Thus it remains to show
$pr_{y}(b^{p}_{v})=\pi_{\beta_{1},\beta_{2}}(b^{p})$ (resp.
$pr_{z}(b^{p}_{v})=\pi_{\tau_{1},\tau_{2}}(b^{p})$) but this is
immediate from the product property of Bernoulli measures and the
definition of the coding map $\pi_{v}$ in section four. \qed The
self similar Bernoulli measures $\pi_{\beta_{1},\beta_{2}}(b^{p})$
resp. $\pi_{\beta_{1},\beta_{2}}(b^{p})$ where extensively studied
in \cite{[NE1]}. We have results on absolut continuity,
singularity and the dimension of this measures.
\section{Measures of full Dimension}
In this section we will ask the question whether there exists
ergodic measures of full dimension for linear solenoids. Our
result is that in general such a measure does not exist, the
dimension of an attractor can not even be approximated by the
dimension of ergodic measures. On the other hand we will proof
that under the assumption of certain symmetries of the system the
equal weighted Bernoulli measure on the attractor has full
dimension. Our theorem is a consequence of both our results on the
dimension of the attractor in section three and our results on the
dimension of ergodic measures in section 6.
\begin{theorem}
For all $v=(\beta_{1},\beta_{2},\tau_{1},\tau_{2})\in (0,1)$ with
$\tau_{1}+\tau_{2}<\beta_{1}+\beta_{2}<1$ we have
\[ \dim b^{0.5}_{v}=\dim_{B}\Lambda_{v}=\dim_{H}\Lambda_{v} \]
if $\beta_{1}=\beta_{2}$ and
\[ \sup\{\dim\mu|\mu\in
\mathcal{E}(\Lambda_{v},f_{v})\}<\dim_{B}\Lambda_{v}=\dim_{H}\Lambda_{v}\]
if
$\beta_{1}\not =\beta_{2}$.\\~\\
For almost all $\beta_{1},\beta_{2}\in (0,0.649)$ with
$\beta_{1}+\beta_{2}\ge 1$ and all $\tau_{1},\tau_{2}\in (0,1)$
with $\tau_{1}+\tau_{2}<\beta_{1}+\beta_{2}$ we have
\[ \dim b^{0.5}_{v}=\dim_{B}\Lambda_{v}=\dim_{H}\Lambda_{v} \]
if $\log_{\tau_{2}}2\beta_{2}=\log_{\tau_{1}}2\beta_{1}$ and
\[ \sup\{\dim\mu|\mu\in
\mathcal{E}(\Lambda_{v},f_{v})\}<\dim_{B}\Lambda_{v}=\dim_{H}\Lambda_{v}\]
if $\log_{\tau_{2}}2\beta_{2}\not =\log_{\tau_{1}}2\beta_{1}$.
\end{theorem}
\proof First assume $\tau_{1}+\tau_{2}<\beta_{1}+\beta_{2}<1$ and
let $d$ be the solution of $\beta_{1}^{d}+\beta_{2}^{d}=1$. In the
case $\beta_{1}=\beta_{2}=1/2$ we have
$\beta_{1}^{d}=\beta_{2}^{2}=1/2$. It is well know in dimension
theory that $\dim\pi_{\beta_{1},\beta_{2}}(b^{0.5})=d$, see for
instance chapter 5 of \cite{[PE]} (the result follows from 15.4 of
this work). Thus we get by Corollary 6.1
\[ \dim
b_{v}^{0.5}=1+\frac{\beta_{1}^{d}\log\beta_{1}^{d}+\beta_{2}^{d}\log\beta_{2}^{d}}{\beta_{1}^{d}\log\tau_{1}+\beta_{2}^{d}\log\tau_{2}}
+(1-\frac{\beta_{1}^{d}\log\beta_{1}+\beta_{2}^{d}\log\beta_{2}}{\beta_{1}^{d}\log
\tau_{1}+\beta_{2}^{d}\log\tau_{2}})d=d+1
\]
But this is by theorem 3.1 the dimension of the attractor
$\Lambda_{v}$. In the case $\beta_{1}\not=\beta_{2}$ we know that
\[ \dim\pi_{\beta_{1},\beta_{2}}(b^{0.5})=\frac{\log 2}{-(0.5\log\beta_{1}+0.5\log\beta_{2})}<d,\]
see again 15.4 of \cite{[PE]}. Consider a sequence $\mu_{n}\in
\mathcal{E}(\Sigma,\sigma)$ with $\mu_{n}\longmapsto b^{0.5}$. By
theorem 6.1 we get
\[ \overline{\lim}_{n\longrightarrow\infty}\dim_{H}(\mu_{n})^{s}_{v}<d+\frac{-2\log2-\log\beta_{1}^{d}-\log\beta_{2}^{d}}{\log\tau_{1}+\log\tau_{2}}<d\]
We thus see that there is a weak$^{\star}$ neighborhood $U$ of
$b^{0.5}$ with
\[ \sup\{\dim\mu_{v}^{s}|\mu\in U \} <d \]
On the other hand it is well known in ergodic theory that
\[ \sup\{ h(\mu)| \mu \in \mathcal{E}(\Sigma,\sigma)\backslash U\}<\log 2 \]
and together we get by theorem 3.1
\[ \sup\{\dim\mu_{v}|\mu\in
\mathcal{E}(\Sigma,\sigma)\}<d+1 \] Our claim now follows from
theorem 3.1 and proposition 4.2.
~\\~\\
Now assume $\tau_{1}+\tau_{2}<\beta_{1}+\beta_{2}$ and
$\beta_{1}+\beta_{2} \ge 1$ and let $d$ be the solution of
$\beta_{1}\tau^{d}+\beta_{2}\tau_{2}^{d}=1$. In the case
$\log_{\tau_{2}}2\beta_{2}=\log_{\tau_{1}}2\beta_{1}$ we have
$\beta_{1}\tau_{1}^{d}=\beta_{2}\tau_{2}^{d}=1/2$. From
\cite{[NE1]} we know that for almost all
$\beta_{1},\beta_{2}\in(0,0.649)$ we have
$\dim\pi_{\beta_{1},\beta_{2}}(b^{0.5})=1$. Thus we get by
Corollary 6.1
\[
\dim_{H}b^{0.5}_{v}
=2+\frac{\beta_{1}\tau_{1}^{d}\log\beta_{1}\tau_{1}^{d}+\beta_{2}\tau_{2}^{d}\log\beta_{2}\tau_{2}^{d}-(\beta_{1}\tau_{1}^{d}\log\beta_{1}+\beta_{2}\tau_{2}^{d}\log\beta_{2})}{\beta_{1}\tau_{1}^{d}\log\tau_{1}+\beta_{2}\tau_{2}^{d}\log\tau_{2}}
\]
\[=2+\frac{\beta_{1}\tau_{1}^{d}\log\tau_{1}^{d}+\beta_{2}\tau_{2}^{d}\log\tau_{2}^{d}}{\beta_{1}\tau_{1}^{d}\log\tau_{1}+\beta_{2}\tau_{2}^{d}\log\tau_{2}}=d+2.
\]
But this is by theorem 3.1 the dimension of the attractor
$\Lambda_{v}$. Now let the
$\log_{\tau_{2}}2\beta_{2}\not=\log_{\tau_{1}}2\beta_{1}$. We get
by theorem 6.1 the following upper estimate
\[ \dim_{H}\mu_{v}^{s}\le 1
-\frac{h_(\mu)+\Xi^{\mu}_{\beta_{1},\beta_{2}}}{\Xi^{\mu}_{\tau_{1},\tau_{2}}}\]
for all $\mu\in \mathcal{E}(\Sigma,\sigma)$. If
$\mu_{n}\longmapsto b^{0.5}$ this yield
\[ \overline{\lim}_{n\longrightarrow\infty}\dim_{H}(\mu_{n})^{s}_{v}\le 1-\frac{\log 2+0.5\log\beta_{1}+0.5\log\beta_{2}}{0.5\log\tau_{1}+0.5\log\tau_{2}}<d+1.\]
We thus again see that there is a weak$^{\star}$ neighborhood $U$
of $b^{0.5}$ with
\[ \sup\{\dim\mu_{v}^{s}|\mu\in U \} <d+1 \]
Now our claim follows in the same way as in the case
$\beta_{1}+\beta_{2}<1$.
 \qed Clearly the statement of theorem 7.1 holds as well if we
interchange the role of $\beta$ and $\tau$. At the end of this
work we again remark that it should be possible to replace the
bound $0.649$ in theorem 7.1 by $1$ using new ideas on continuity
of self-similar Bernoulli measures
$\pi_{\beta_{1},\beta_{2}}(b^{0.5})$ in the case
$\beta_{1}+\beta_{2}\ge 1$.
\section{Appendix: General facts in dimension theory} We will
here first define the most important quantities in dimension
theory and collect some basic facts we need. We refer to the book
of Falconer \cite{[FA]} and the book of Pesin \cite{[PE]} for a
more detailed discussion of dimension theory.\newline\newline Let
$Z\subseteq\mathbb{R}^{q}$. We define the {\it $s$-dimensional
Hausdorff measure} $H^{s}(Z)$ of $Z$ by
\[ H^{s}(Z)=\lim_{\lambda \longrightarrow 0}\inf\{\sum_{i\in I}(\mbox{diam} U_{i})^{s}|Z\subseteq\bigcup_{i\in I}U_{i}\mbox{ and }\mbox{diam}(U_{i})\le\lambda\}. \]
The {\it Hausdorff dimension} $\dim_{H}Z$ of $Z$ is given by
\[ \dim_{H}Z=\sup\{s|H^{s}(Z)=\infty\}=\inf\{s|H^{s}(Z)=0\}.\]
Let $N_{\epsilon}(Z)$ be the minimal number of balls of radius
$\epsilon$ that are needed to cover $Z$. We define the {\it upper
box-counting dimension} $\overline{\dim}_{B}$ resp. {\it lower
box-counting dimension} $\overline{\dim}_{B}$ of $Z$ by
\[ \overline{\dim}_{B}Z=\overline{\lim}_{\epsilon\longrightarrow 0} \frac{\log N_{\epsilon}(Z)}{-\log\epsilon}\qquad
\underline{\dim}_{B}Z=\underline{\lim}_{\epsilon\longrightarrow 0}
\frac{\log N_{\epsilon}(Z)}{-\log\epsilon}.\]
\begin{proposition} If $Z\subseteq\mathbb{R}^{q}$ and $I\subseteq\mathbb{R}$ is an interval
then $\dim_{H}(Z\times I)=\dim_{H}+1$ and $\dim_{B}(Z\times
I)=\dim_{B}+1$ holds for both upper an lower dimension.
\end{proposition}
The statement for the Hausdorff dimension follows from proposition
7.4. of \cite{[FA]} and the statement for the box-counting
dimension is easy to see using 3.1. of
\cite{[FA]}.\newline\newline Now let $\mu$ be a Borel probability
measure on $\mathbb{R}^{q}$. We define the dimensional theoretical
quantities for $\mu$ by
\[ \dim_{H}\mu=\inf\{\dim_{H}Z|\mu(Z)=1\}\]
and
\[ \overline{\underline{\dim}}_{B}\mu=\lim_{\rho\longrightarrow 0}
\inf\{\overline{\underline{\dim}}_{B}Z|\mu(Z)\ge 1-\rho \}.\] We
introduce one more notion of dimension for a measure $\mu$.  The
{\it upper local dimension} $\overline{d}(x,\mu)$ resp. {\it lower
local dimension} $\underline{d}(x,\mu)$ of the measure $\mu$ in a
point $x$ is defined by
\[ \overline{d}(x,\mu)=\overline{\lim}_{\epsilon\longrightarrow 0} \frac{\mu(B_{\epsilon}(x)) }{\log\epsilon}\qquad
\underline{d}(x,\mu)=\underline{\lim}_{\epsilon\longrightarrow 0}
\frac{\mu(B_{\epsilon}(x))}{\log\epsilon}.\] The relations between
the local dimension and the other notion of dimension of measures
are described in the following theorem:
\begin{theorem}
Let $\mu$ be a Borel probability measure on $\mathbb{R}^{q}$. If
\[ \overline{d}(x,\mu)=\underline{d}(x,\mu)=c\] almost for holds
for almost all $x\in\mathbb{R}^{q}$ we have
\[ dim_{H}\mu=\dim_{B}\mu=c.\]
\end{theorem}
A proof of this theorem is contained in the work of Young
\cite{[YO]}. If the conditional in this theorem holds, the measure
$\mu$ is called {\it exact dimensional} and the common value of
the dimensions is denoted by $\dim\mu$ and maybe called the
fractal dimension of the measure. 
\end{document}